\documentclass[onecolumn,final,10pt, reqno]{amsbook}
\pdfoutput=1
\usepackage{enumerate} 
\usepackage{graphicx} 
\usepackage{lscape} 
\usepackage{amssymb}

\usepackage{hyperref} 

\makeatletter
\def\specialsection{\@startsection{section}{1}%
  \z@{\linespacing\@plus\linespacing}{.5\linespacing}%
  {\normalfont}}
\def\section{\@startsection{section}{1}%
  \z@{.7\linespacing\@plus\linespacing}{.5\linespacing}%
  {\large\scshape\bfseries}}
\makeatother

\usepackage{etoolbox}
\makeatletter
\patchcmd{\@maketitle}{\newpage}{}{}{} 
\makeatother

\makeatletter
\g@addto@macro{\maketitle}{\clearpage
\begin{minipage}[h]{0.9\textwidth}
\begin{bfseries}
\begin{center}
\Large
	Copyright
\normalsize
\end{center}
\end{bfseries}
\vspace{0.25in}
The underlying manuscript is held by the National Archives in the UK  and can be accessed at \url{www.nationalarchives.gov.uk} using reference number HW 25/37. Readers are encouraged to obtain a copy. 
\\

The original  work was under Crown copyright, which has now expired; the work is now in the public domain. 
\end{minipage}
}
\makeatother

\title{The Applications of Probability to Cryptography}
\author{Alan M. Turing}

\begin{document}
\frontmatter
\maketitle
\begin{bfseries}
\begin{center}
\Large
	Editor's Notes
\normalsize
\end{center}
\begin{flushleft}
\textsc{Provenance}
\end{flushleft}
\end{bfseries}

Two Second World War research papers by Alan Turing were declassified recently.  The papers, \textit{The Applications of Probability to Cryptography} and its shorter companion \textit{Paper on Statistics of Repetitions}, are available from from the National Archives in the UK at \url{www.nationalarchives.gov.uk.} 

The released papers give the full text, along with figures and tables, and provide a fascinating insight into the preparation of the manuscripts, as well as the style of writing at a time when typographical errors were corrected by hand, and mathematical expression handwritten into spaces left in the text.

Working with the papers in their original format provides some challenges, so they have been typeset for easier reading and access. We recommend that the typeset versions are read with a copy of the original manuscript at hand.

This document contains the text and figures for \textit{The Applications of Probability to Cryptography}, the companion paper is also available in typeset form from arXiv at \url{www.arxiv.org/abs/1505.04715}. These notes apply to both documents.

Separately, a journal article by Zabell\footnote{Zabell, S. 2012. ``Commentary on Alan M.Turing: The Applications of Probability to Cryptography''  \textit{Cryptologia}, 36:191-214.} provides an analysis of the papers and further background information.
\smallskip
\begin{bfseries}
\begin{flushleft}
\textsc{The text}
\end{flushleft}
\end{bfseries}

It is not our intent to cast Alan Turing's manuscripts into a journal style article, but more to provide clearer access to his writing and, perhaps, to answer the questions ``If Turing had have had access to typesetting software, what would his papers have looked like?''. Consequently no ``house-style''  copy-editing has been imposed. Occasional punctuation has been added to improve readability, some obvious errors corrected, and paragraph breaks added to ease the reading of long text blocks - and occasionally to give a better text flow. Turing uses typewriter underlining, single, and double quotes to indicate emphasis or style; these have been implemented using font format changes, double quotes are used as needed. 

The manuscript has many typographical errors, deletions, and substitutions, all of which are indicated by over-typing, crossed out items, and handwritten pencil or ink annotations. These corrections have been implemented in this document to give the text that we presume Turing intended. Additionally, there are some hand written notes in the manuscript, which may or may not be by Turing; these are indicated by the use of footnotes.

British English spelling is used in the manuscript and this is retained, so words such as favour, neighbourhood, cancelling, etc. will be encountered. Turing appears to favour the spellings bigramme, trigramme, tetragramme, etc., although he is not always consistent; throughout this document the favoured rendering is used.

Turing's wording is unchanged to give the full flavour of his original style. This means that ``That is to say that we suppose for instance that ..... '' will be encountered - amongst others!

Both papers end abruptly, no summary or conclusion is offered, perhaps the papers are incomplete or pages are missing. To indicate the end of the manuscript we have marked the end of each paper with a printing sign - an infinity symbol between two horizontal bars.

In the section on a letter subtractor problem, reference is made to other methods to be discussed later in the paper. This does not happen - perhaps another indicator of an incomplete paper or missing pages.

Finally, Turing uses some forward page references that appear in the manuscript as \textit{see(p  )}, obviously intending to return and complete the reference. This also does not happen, so these references remain unresolved.

In short, we strive to represent Turing's text as he wrote it.

\smallskip
\begin{bfseries}
\begin{flushleft}
\textsc{Ciphertext, cleartext, etc.}
\end{flushleft}
\end{bfseries}

In an attempt to capture the flavour of the time, ciphertext, cleartext, keys, etc. are displayed in a fixed pitch, bold, non-serif font to represent the typewriter, teletype, and telegraph machines that would have printed the original code, \textit{viz.}~\texttt{CONDITIONS}.
 
\smallskip
\begin{bfseries}
\begin{flushleft}
\textsc{Mathematics}
\end{flushleft}
\end{bfseries}

In the manuscript all mathematics is hand written in ink and pencil in spaces left between the typed text. Sometimes adequate space was left, other time not, and the handwriting spills into margins and adjacent lines, adding to the reading challenge. We have cast all mathematics into standard in-line or display formats as appropriate. We have used the \texttt{mathcal }font in places to capture the flavour of Turing's handwriting, \textit{e.g.} ``the probability \textit{p}'' appears as ``the probability $\mathcal{P}$''.

Turing uses no punctuation in his mathematics, this has been added to be consistent with modern practice\footnote{See, for instance, Higham, Nicholas J. 1998. ``Handbook of writing for the Mathematical Sciences'', SIAM, Philadelphia.}; he also uses letters to reference equations -  numbers are used in this document. In many places we have added parentheses to give clarity to an expression, and in some places where Turing is inconsistent in his uses of parentheses for a mathematical phrase (the expression for letter probability in the Vigen\`{e}re in particular) we have chosen one format and been consistent in its use.

As Turing demonstrates a love of dense mathematics the algebraic multiplication symbol $\times$ has occasionally been used for readability, so all standard forms of multiplication will be encountered, \textit{viz.,} $ab, a \times b, a \cdot b$. Finally, convention suggests that the subject of a formula or expression sits on its own on the left hand side of the equals sign, with the subsidiary variables collected on the right hand side. Turing adheres to this convention as it suits him, his preference is retained.
 
In short, we strive to retain the elegance of Turing's mathematics, whilst casting it into a modern format.

\smallskip
\begin{bfseries}
\begin{flushleft}
\textsc{Figures and tables}
\end{flushleft}
\end{bfseries}

All figures have been included with rearrangement of some items to improve clarity or document flow. Turing uses a variety of papers, styles, inks, pen, and pencil; these have all been represented in standard figure and table format.

\smallskip
\begin{bfseries}
\begin{flushleft}
\textsc{Contents page}
\end{flushleft}
\end{bfseries}

Turing provides a rudimentary Contents for \textit{The Applications of Probability to Cryptography}, this has been reworked with some additions to make it more meaningful. \textit{Paper on Statistics of Repetitions}, being much shorter, requires no~Contents.

\smallskip
\begin{bfseries}
\begin{flushleft}
\textsc{Editors}
\end{flushleft}
\end{bfseries}

The editor can be contacted at: \textit{ian.taylor@maths.oxon.org}.
\tableofcontents
\mainmatter
%
%
\chapter{Introduction}
\section{Preamble}
The theory of probability may be used in cryptography with most effect when the type of cipher used is already fully understood, and it only remains to find the actual keys. It is of rather less value when one is trying to diagnose the type of cipher, but if definite rival theories about the type of cipher are suggested it may be used to decide between them.

\section{Meaning of probability and odds}
I shall not attempt to give a systematic account of the theory of probability, but it may be worth while to define shortly \textit{probability} and \textit{odds}. The \textit{probability} of an event on certain evidence is the proportion of cases in which that event may be expected to happen given that evidence. For instance if it is known the 20\% of men live to the age of 70, then knowing of Hitler only \textit{Hitler is a man} we can say that the probability of Hitler living to the age of 70 is 0.2. Suppose that we know that \textit{Hitler is now of age 52} the probability will be quite different, say 0.5, because 50\% of men of 52 live to 70.

The \textit{odds} of an event happening is the ratio $\mathcal{P}/(1 - \mathcal{P})$ where $\mathcal{P}$ is the probability of it happening. This terminology is connected with the common phraseology \textit{odds of 5:2 on} meaning in our terminology that the  odds are 5/2.

%
%
\section{Probabilities based on part of the evidence}
When the whole evidence about some event is taken into account it may be extremely difficult to estimate the probability of the event, even very approximately, and it may be better to form an estimate based on a part of the evidence, so that the probability may be more easily calculated. This happens in cryptography in a very obvious way. The whole evidence when we are trying to solve a cipher is the complete traffic, and the events in question are the different possible keys, and functions of the keys. Unless the traffic is very small indeed the theoretical answer to the problem ``\,What are the probabilities of the various keys?\,'' will be of the form ``\,The key $\dots$ has a probability differing almost imperceptibly from 1 (certainty) and the other keys are virtually impossible''. But a direct attempt to determine these probabilities would obviously not be a practical method.

\section{\textit{A priori} probabilities}
The evidence concerning the possibility of an event occurring usually divides into a part about which statistics are available, or some mathematical method can be applied, and a less definite part about which one can only use one's judgement. Suppose for example that a new kind of traffic has turned up and that only three messages are available. Each message has the letter \texttt{V} in the 17th place and \texttt{G} in the 18th place. We want to know the probability that it is a general rule that we should find \texttt{V} and \texttt{G} in these places. We first have to decide how probable it is that a cipher would have such a rule, and as regards this one can probably only guess, and my guess would be about $1/5,000,000$. 
%
%
This judgement is not entirely a guess; some rather insecure mathematical reasoning has gone into it, something like this:-

The chance of there being a rule that two consecutive letters somewhere after the 10th should have certain fixed values seems to be about $1/500$ (this is a complete guess). The chance of the letters being the 17th and 18th is about $1/15$ (another guess, but not quite as much in the air). The probability of a letter being \texttt{V} or \texttt{G} is $1/676$ (hardly a guess at all, but expressing a judgement that there is no special virtue in the bigramme \texttt{VG}). Hence the chance is $1/(500 \times 15 \times 676)$ or about $1/5,000,000$. This is however all so vague, that it is more usual to make the judgment ``$1/5,000,000$'' without explanation.

The question as to what is the chance of having a rule of this kind might of course be resolved by statistics of some kind, but there is no point in having this very accurate, and of course the experience of the cryptographer itself forms a kind of statistics.

The remainder of the problem is then solved quite mathematically. Let us consider a large number of ciphers \textit{chosen at random}. $N$ of them say. Of these $N/5,000,000$ of them will have the rule in question, and the remainder not. Now if we had three messages of each of the ciphers before us, we should find that for  each of the ciphers with the rule, three messages have \texttt{VG} in the required place, but of the remaining $(4,999,999 \times N) / 5,000,000$ only a proportion $1/676^3$ will have them. Rejecting the ciphers which have not the required characteristics we are left with 
%
%
$N/5,000,000$ cases where the rule holds, and  $\left( 4,999,999 \times N \right) /( 5,000,000 \times 676^3)$ cases where it does not. This selection of ciphers is a random selection of ones which have all the known characteristics of the one in question, and therefore the odds in favour of the rule holding are:
\begin{gather*}
	 \frac{N}{5,000,000} : \frac{4,999,999 \times N} {5,000,000 \times 676^3}, \\
	   i.e \quad 676^3:4,999,9999, \\
	   \text{or about} \  60:1 \  \text{on}. 
\end{gather*}
It should be noticed that the whole argument is to some extent fallacious, as it is assumed that there are only two possibilities, \textit{viz.} that either \texttt{VG} must always occur in that position, or else that the letters in the 17th and 18th positions are wholly random. There are however many other possibilities worth consideration, \textit{e.g.}
\begin{enumerate}
	\item On the day in question we have \texttt{VG} in the position in question. 
	\item Or on another day we have some other fixed pair of letters. 
	\item Or in the positions 17, 18 we have to have one of the four combinations \texttt{VG}, \texttt{RH}, \texttt{OM}, \texttt{IL} and by chance \texttt{VG} has been chosen for all the three messages we have had. 
	\item Or the cipher is a simple substitution and \texttt{VG} is the substitute of some common bigramme, say \texttt{TH}.	
\end{enumerate}

The possibilities are of course endless, and it is therefore always necessary to bear in mind the possibility of there being other theories not yet suggested.

The \textit{a priori} probability sometimes has to be estimated as above by some sort of guesswork, but often the situation is more satisfactory. Suppose for example that we know that a certain cipher is a simple substitution, the keys having no specially noticeable properties. Suppose also that we have 50 letters of such a message including five occurrences of \texttt{P}. We want to know how probable it it that \texttt{P} is the substitute of \texttt{E}. As before we have to answer two questions. 
%
%
\begin{enumerate}
	\item How likely is it that \texttt{P} would be the substitute of \texttt{E} neglecting the evidence of the five \texttt{E}s occurring in the message?
	\item How likely are we to get 5 \texttt{P}s?
	\begin{enumerate}
		\item If \texttt{P} is not the substitute of \texttt{E}
		\item If \texttt{P} is the substitute of \texttt{E}.
	\end{enumerate}
\end{enumerate}
I will not attempt to answer the second question for the present. The answer to the first is simply that the probability of a letter being the substitute of \texttt{E} is independent of what the letter is, and is therefore always $1/26$, in particular it is $1/26$ for the letter \texttt{P}. The only guesswork here is the judgement that the keys are chosen at random.

\section{The Factor Principle}
Nearly all applications of probability to cryptography depend on the \textit{factor principle} (or Bayes' Theorem). This principle may first be illustrated by a simple example. Suppose that one man in five dies of heart failure, and that of the men who die of heart failure two in three die in their beds, but of the men who die from other causes only one in four dies in their beds. (My facts are no doubt hopelessly inaccurate). Now suppose we know that a certain  man died in his bed. What is the probability that he died of heart failure? Of all numbering \textit{N} say we find that
\begin{align*}
	N \times &(1/5) \times (2/3)  &&\textit{ die in their beds of heart failure } \\
	N \times &(1/5 \times (1/3)  &&\dots \textit{  elsewhere } \quad \dots \dots \dots \dots \dots \dots \\
	N \times &(4/5) \times (1/4)  &&\textit{ die in their beds from other causes } \\
	N \times &(4/5 \times (3/4)  &&\dots \textit{  elsewhere } \quad \dots \dots \dots \dots \dots  \dots
\end{align*}
Now as our man died in his bed we do not need to consider the cases of men who did not die in their beds, and these consist of
\begin{align*}
	N \times &(1/5) \times (2/3)  &&\textit{ cases of heart failure and  } \\
	N \times &(4/5) \times (1/4)  &&\textit{ from other causes }
\end{align*}
and therefore the odds are $1 \times  (2/3) : 4 \times  (1/4) $in favour of heart failure. If this had been done algebraically the result would have been
%
%
\begin{multline*}
	\textit{A posteriori odds of the theory} \\
	= \textit{A priori odds of the theory} \\
	\times \frac{\textit{Probability of the data being fulfilled if the theory is true}}{\textit{Probability of the data being fulfilled if the theory is false}}.
\end{multline*}
In this the \textit{theory} is that the man died of heart failure, and the \textit{data} is that he died in his bed. 

The general formula above will be described as the \textit{factor principle}, the ratio
\[
	\frac{\textit{Probability of the data if the theory is true}}{\textit{Probability of the data if the theory is false}},
\]
\\
is called the factor for the theory on account of the data.

\section{Decibanage}
Usually when we are estimating the probability of a theory there will be several independent pieces of evidence \textit{e.g.} following our last example, where we want to know whether a certain man died of heart failure or not, we may know
\begin{enumerate}
	\item He died in his bed
	\item His father died of heart failure
	\item His bedroom was on the ground floor
\end{enumerate}
and also have statistics telling us
\begin{enumerate}[(a)]
	\item 2/3 of men who die of heart failure die in their beds
	\item 2/5 \dots \dots \dots \dots \dots \dots \dots \dots \dots \dots \dots have fathers who died of heart failure
	\item 1/2 \dots \dots \dots \dots \dots \dots \dots \dots \dots \dots \dots have bedroom on the ground floor
	\item 1/4 of men who died from other causes die in their beds
	\item 1/6 \dots \dots \dots \dots \dots \dots \dots \dots \dots \dots \dots have fathers who died of heart failure
	\item 1/20 of men who die of other cause have their bedrooms on the ground floor
\end{enumerate}

%
%
Let us suppose that the three pieces of evidence are independent of one another if we know that he died of heart failure, and also if we know that he did not die of heart failure. That is to say that we suppose for instance that knowing that he slept on the ground floor does not make it any more likely that he died in his bed if we knew all along that he died of heart failure. When we make these assumptions the probability of a man who died of heart failure satisfying all three conditions is obtained simply by multiplication, and is $(2/3) \times (2/5) \times (1/2)$ and likewise for those who died from other causes the probability is $(1/4) \times (1/6) \times (1/20)$, and the factor in favour of the heart theory failure is
\[
	\frac{(2/3) \times (2/5) \times (1/2)}{(1/4) \times (1/6) \times (1/20)}.
\]
We may regard this as the product of three factors $(2/3)/(1/4)$ and $(2/5)/(1/6)$ and $(1/2)/(1/20)$ arising from from the three independent pieces of evidence. Products like this arise very frequently, and sometimes one will get products involving thousands of factors, and large groups of these factors may be equal. We naturally therefore work in terms of the logarithms of the factors. The logarithm of the factor, taken to the base $10^{1/10}$ is called \textit{decibanage in favour of the theory}. A \textit{deciban} is a unit of evidence; a piece of evidence is worth a deciban if it increase the odds of the theory in the ratio $10^{1/10} : 1$. The deciban is used as a more convenient unit that the \textit{ban}. The terminology was introduced in honor of the famous town of Banbury. 

%
%
Using this terminology we might say that the fact that our man died in bed scores 4.3 decibans in favour of the heart failure theory $(10 \log (8/3) = 4.3)$. We score a further 3.8 decibans for his father dying of heart failure, and 10 for his having his bedroom on the ground floor, totalling 18.1 decibans. We then bring in the \textit{a priori} odds 1/4 or $10^{-6/10}$ and the result is the the odds are $10^{12.1/10}$, or as we may say ``12.1 deciban up on evens''. This means about 16:1 on.

%
%
\chapter{Straightforward Cryptographic Problems}
\section{Vigen\`{e}re}
The factor principle can be applied to the solutions of a Vigen\`{e}re problem with great effect. I will assume here that the period of the cipher has already been determined. Probability theory may be applied to this part of the problem also, but that is not so elementary. Suppose our cipher, written out in its correct period~is\footnote{ Turing's statement of the ciphertext is slightly different to what he decodes. The N M at the end the first line are reversed to read \texttt{DKQHSHZMNP} in Fig 5, which gives the correct cleartext.}
\vspace{0.1in}
\Large
\begin{bfseries}
	\begin{ttfamily}
		\begin{center}
			D K Q H S H Z N M P  \\
			R C V X U H T E A Q \\
			X H P U E P P S B K \\
			T W U J A G D Y O J  \\
			T H W C Y D Z H G A  \\
			P Z K O  X O E Y A E  \\
			B O K B U B P I K R  \\
			W W A C E J P H L P  \\
			T U Z Y F H L R Y C
		\end{center}
	\end{ttfamily}
\end{bfseries}
\normalsize
\begin{rmfamily}
	\begin{center}
		\textsc{Figure 1.} Vigen\`{e}re problem. \\
		\textit{(It is only by chance that it makes a rectangular array.)}
	\end{center}
\end{rmfamily}
\vspace{4pt}
Let us try to find the key for the first column, and for the moment let us only take into account the evidence afforded by the first letter \texttt{D}. Let us first consider the key \texttt{B}. The factor principle tells us
\begin{multline*}
	\textit{Odds in favour of key \texttt{B}} = \textit{ A priori odds in favour of key  \texttt{B}} \\
	\times \frac{\textit{Probability of getting \texttt{D} in cipher if key is \texttt{B}}}{\textit{Probability of getting key \texttt{D} in cipher if key is not \texttt{B}}}
\end{multline*}

Now the \textit{a priori} odds in favour of key \texttt{B} may be taken as 1/25. The probability of getting \texttt{D} in the cipher with the key \texttt{B} is just the probability of getting \texttt{C} in the clear which (using the count on 1000 letters in Fig 2) is 0.021. If however the key is not \texttt{B} we can have any letter other the \texttt{C} in the clear, and the probability is \linebreak (1 - 0.021)/25. Using the evidence of the \texttt{D} then the odds in favour of the key \texttt{B} are
\[
	\frac{1}{25} \times \left (\frac{25 \times 0.021}{1 - 0.021} \right).
\]
%
%
We may then consider the effect of the next letter in the column \texttt{R} which gives a further factor of (25 x 0.064)/(1 - 0.064). We are here assuming that the evidence of the \texttt{R} is independent of the evidence of the \texttt{D}. This is not quite correct, but is a useful approximation; a more accurate method of calculation will be given later. Let us write $\mathcal{P}_{\alpha}$ for the frequency of the letter $\alpha$ in plain language. Then our final estimate for the odds in favour of key \texttt{B} is
\[
	\frac{1}{25} \prod_{i} \frac{25 \mathcal{P}_{\alpha_{i} - 1}}{1 - \mathcal{P}_{\alpha_{i} - 1}} .
\]
where $\alpha_{1}, \alpha_{2},  \dots$ is the series of letters in the 1st column, and we use the letters and numbers interchangeably, \texttt{A} meaning 1, \texttt{B} meaning 2, $\dots$, \texttt{Z} meaning 26 or 0. More generally for key $\beta$ the odds are
\[
	\frac{1}{25}  \prod_{i} \frac{25\mathcal{P} _{ \alpha_{i} - \beta + 1}}{1 - \mathcal{P}_{ \alpha_{i} - \beta + 1}}.
\]
The value of this can be calculated by having a table of the decibanage corresponding to the factors $25 \mathcal{P}_{\alpha} / (1 - \mathcal{P}_{\alpha})$. One then decodes the column with the various possible keys, looks up the decibanage, and adds them up.

The most convenient form for doing this is a table of values of $20 \log_{10}[25 \mathcal{P}_{\alpha} / (1 - \mathcal{P}_{\alpha})]$, taken to the nearest integer, or as we may say, the values of the score in  \textit{half decibans}. One may also have columns showing multiples of these, and the table made of double height\footnote{ Turing provides a table of double height for Fig 3 to allow the ``gadget'' of Figure 4 to be used with any letter of the alphabet as a decode key - hence the double alphabet. Figure 4 can be prepared as a transparency, with the original markings cleared, and markings for the new decode letter added. Fig 3 and Fig 4 are correctly proportioned in this document for this to work.} (Fig \ref{Fi:FIG3}). For the first column with key \texttt{B} the decoded column is \texttt{CQWS\textbullet \textbullet OAV},\footnote{ \texttt{S\textbullet\textbullet} means \texttt{SSS}, for a total of three letter \texttt{S}, as noted in the following arithmetic. The linear decode for the example is \texttt{CQWSSOAVS}} and we score -5 for \texttt{C}, -26 for \texttt{Q}, -5 for \texttt{W}, 17 for the three letters \texttt{S}, 5 for \texttt{O}, 7 for \texttt{A} and -10 \texttt{V}, totalling -17. These calculations can be done very quickly by the use of the transparent gadget Fig \ref{Fi:FIG4} , in which squares are ringed in pencil to show the number of letters occurring in the column.
\vspace{8pt}
%
%
\begin{center}
	\begin{tabular}{lrlrlr}
		A	& 84 &J	&2	&S	& 73	  \\
		B	& 23 &K	&5	&T	& 81	\\
		C	& 21	 &L	&38	&U	& 19  \\
		D	& 46	 &M	&34	&V	& 11   \\
		E	&116  &N	&66	&W	& 21  \\
		F	& 20	 &O	&66	&X	& 16  \\
		G	& 25	 &P	&15	&Y	& 24 \\
		H	& 49	 &Q	&2	&Z	& 3  \\
		I	& 76	 &R	&64  
	\end{tabular}
\end{center}
\begin{rmfamily}
	\begin{center}
		\textsc{Figure 2.} Count on 1000 letters. \\
		\textit{(English text)\\
		The value for X has been taken more of less at random as a compromise \\
		between real language \& telegraphese. Also I added to each entry (see p  )\footnote{ Forward reference left unresolved in the manuscript.}.}
	\end{center}
\end{rmfamily}
\setcounter{figure}{2}
\newpage
%
%
\begin{figure}[h]
\centering\includegraphics[scale=0.9]{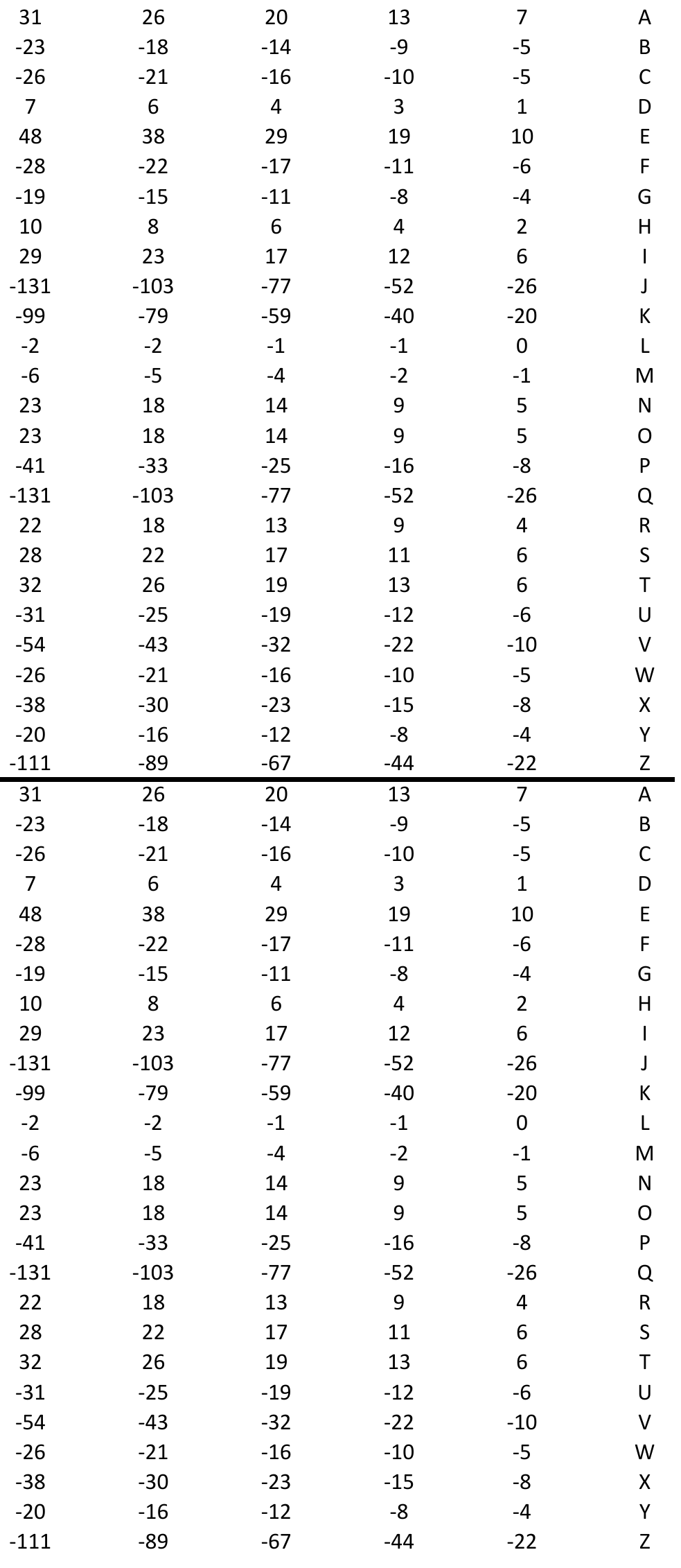}
\caption{Table for scoring a Vigen\`{e}re.}\label{Fi:FIG3}
\begin{center}
		In units of half a deciban. 
\end{center}
\end{figure}
%
%
\newpage
\begin{landscape}
\begin{figure}[h]
\centering\includegraphics[scale=0.85]{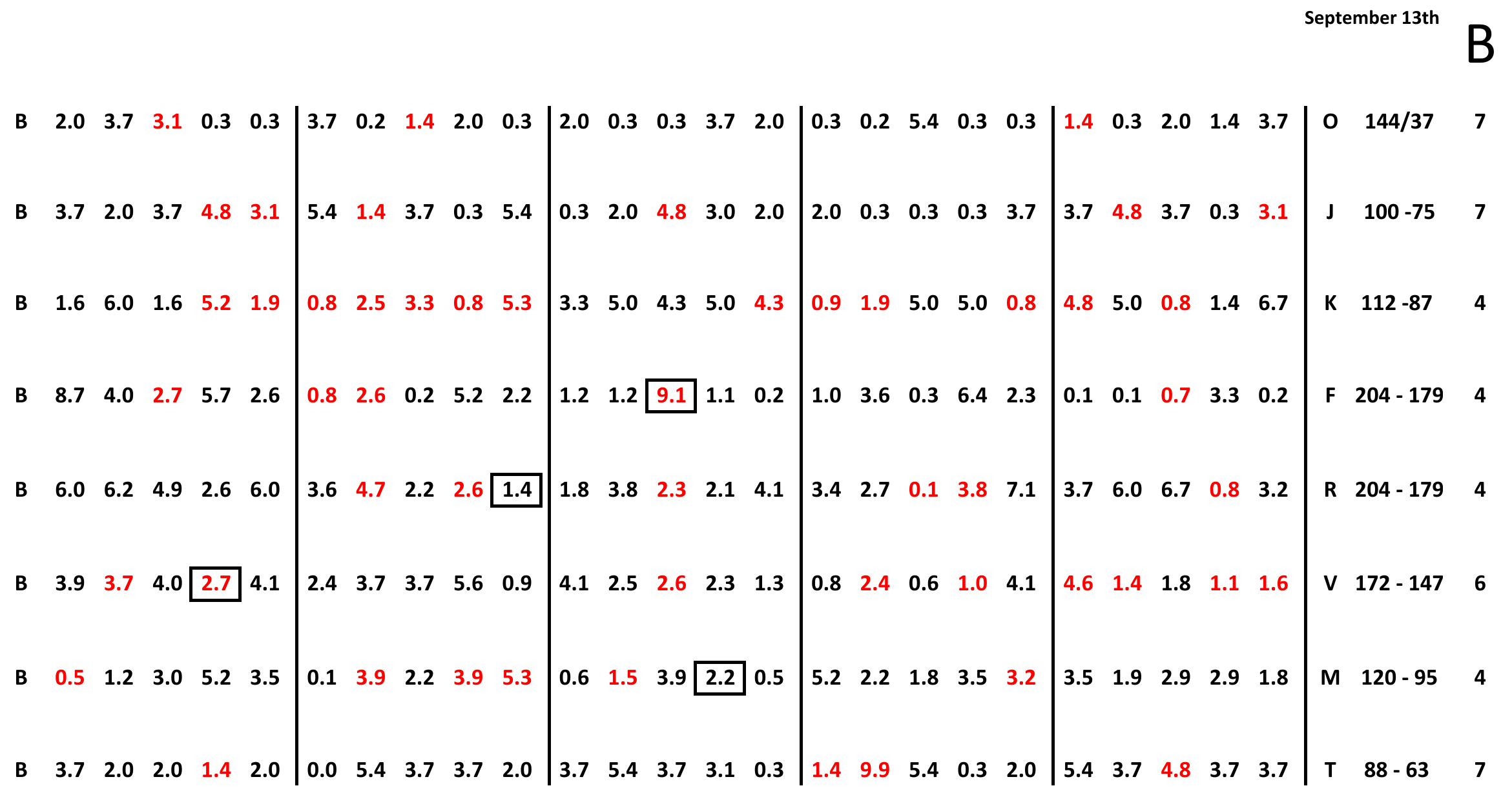}
\begin{center}
\textit{Editor -This table is not referenced in the manuscript, does not have a page number, and \\ nothing appear to indicate its purpose. It is presented here for completeness.}
\end{center}
\end{figure}
\end{landscape}
%
%
\newpage
\begin{figure}[h]
\centering\includegraphics[scale=0.62]{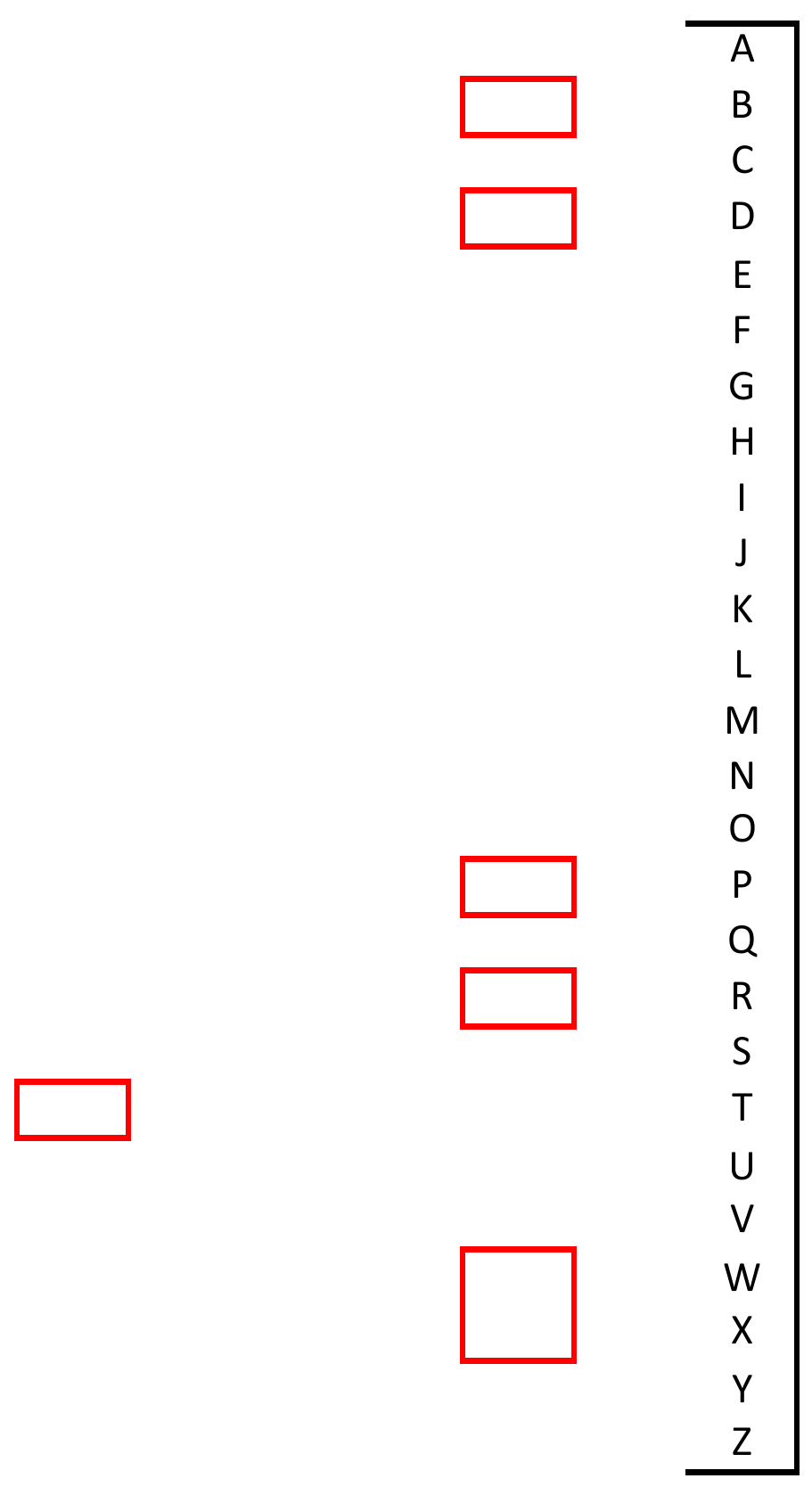}
\caption{Apparatus for scoring a Vigen\`{e}re.}\label{Fi:FIG4}
	\begin{center}
		Pencil marks arranged for 1st wheel of Fig. 1.
	\end{center}
\end{figure}
%
%
The gadget may be placed over Fig \ref{Fi:FIG3} in various positions corresponding to the various keys. The score is obtained by adding up the numbers showing through the various squares. In Fig \ref{Fi:FIG5} the alphabet has been written in a vertical below the cipher text of Fig 1, each letter representing a possible key. The score for each key has been written opposite the key, and under the relevant column. An \texttt{X} denotes a bad score, not worth adding up. Usually these will be -15 or worse. It will be seen that for the first column \texttt{P}, having a score of 43 is extremely likely to be right, especially as there is no other score better than 8. If we neglect this latter fact the odds for the key are $(1/25) 10^{2.15}$ \textit{i.e.} about 5:1 on. The effect of decoding this column with key \texttt{P} has been shown underneath. 

For the second column the best key is \texttt{O}, but is by no means so certain as the first column. The decode for this column is also shown, and provides very satisfactory combinations with the first column, confirming both keys. (This confirmation could also be based on probability theory, given a table of bigramme frequencies). In the third column \texttt{I} and \texttt{C} are best although \texttt{D} would be very possible, and in the fourth column \texttt{Q} and \texttt{U} are best.
%
%
\begin{figure}[hbt]
\centering\includegraphics[scale=0.8]{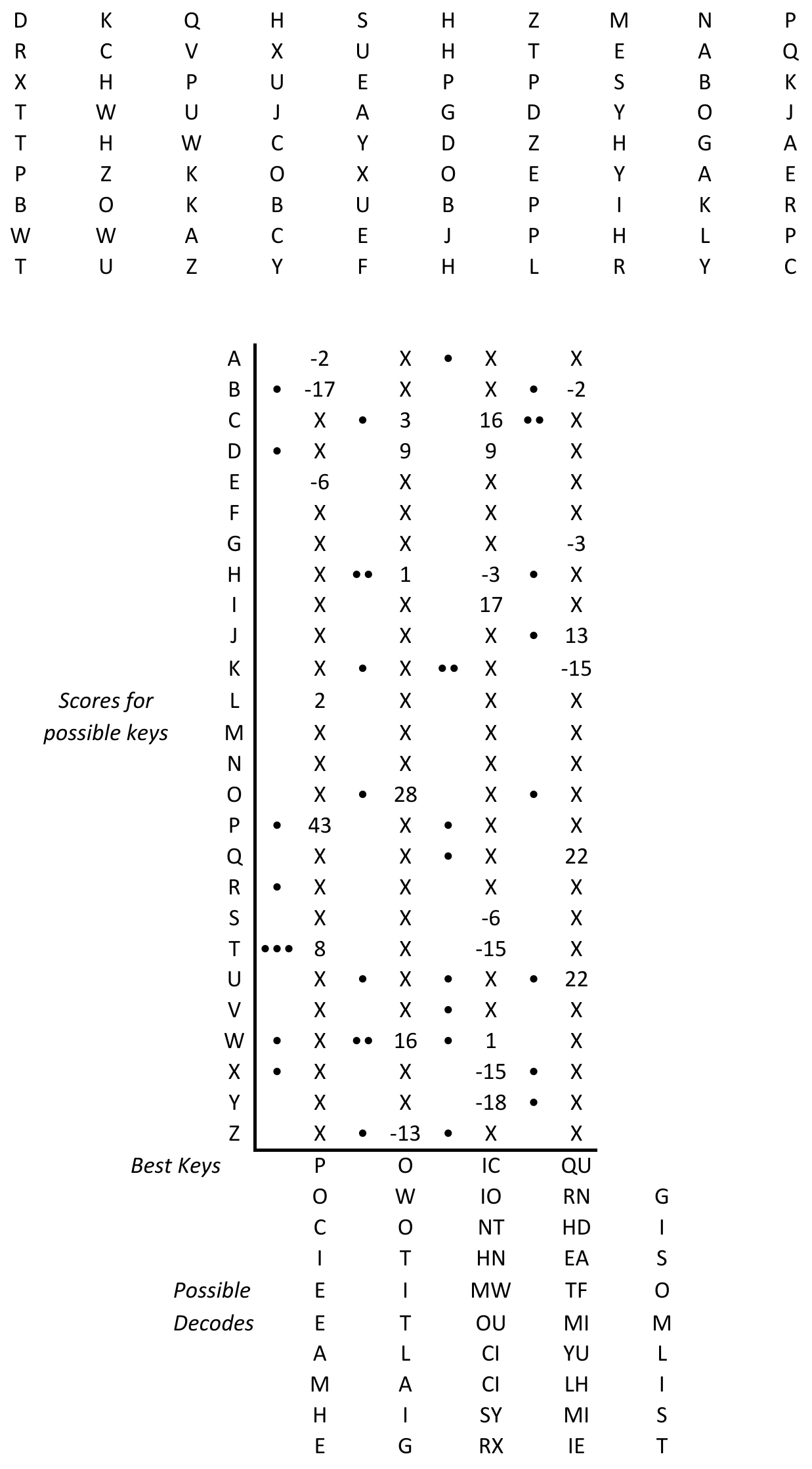}
\caption{Scoring and solving a Vigen\`{e}re.}\label{Fi:FIG5}
\end{figure}
\newpage
%
%
Writing down the possible decodes we see that the first line must read \texttt{OWING} and this makes the other lines read \texttt{CONDI}, \texttt{ITHAS}, \texttt{EIMPO}, \texttt{ETOIM}, \texttt{ALCUL}, \texttt{MACHI}, \texttt{HISIS}, \texttt{EGRET}. By filling in the word \texttt{CONDITIONS} the whole can now be decoded.\footnote{ Solution: Keylength - 10, Key - \texttt{POIUMOLQNY}, Cleartext - \texttt{OWINGTOWAR CONDITIONS ITHASBECOM EIMPOSSIBL ETOIMPORTC ALCULATING MACHINESXT HISISVERYR EGRETTTABLE}}

A more accurate argument would run as follows. For the first column, instead of setting up as rival theories the two possibilities that \texttt{B} is the key and that \texttt{B} is not we can set up 26 rival theories that the key is \texttt{A} or \texttt{B} or \dots  \ \texttt{Z}, and we may apply the factor principle in the form:-
\small
\begin{multline*}
	\frac{\textit{A posteriori probability of key A}}{\textit{A priori probability of key \texttt{A}} \times \textit{Probability of getting the given column with key \texttt{A}}}, \\
	= \frac{\textit{A posteriori probability of key \texttt{B}}}{\textit{A priori probability of key \texttt{B}} \times \textit{Probability of getting the given column with key \texttt{B}}},
\end{multline*}
\begin{flushleft}
	\quad \textit{= etc.}
\end{flushleft}
\normalsize
The argument	to justify this form of factor principle is really the same as for the original form. Let $q_{\beta}$ be the \textit{a priori} probability of key $\beta$. Then out of \textit{N} cases we have $Nq_{\beta}$ cases of key $\beta$. Let $\mathcal{P} \left(\beta, C \right)$ be the probability of getting the column \textit{C} with key $\beta$, then we have rejected the cases where we get columns other than \textit{C} we find that there are $Nq_{\beta} \mathcal{P} \left( \beta, C \right)$ cases of key $\beta$ \textit{i.e.} the \textit{a posteriori} probability of key $\beta$ is $\mathit{K} Nq_{\beta} \mathcal{P} \left( \beta, C \right)$, where $\mathit{K}$ is independent of $\beta$.

We have therefore to calculate the probability of getting the column \textit{C} with key $\beta$ and this is simply $\prod_{i} \mathcal{P}_{\left( \alpha_{i} - \beta + 1 \right)}$, \textit{i.e.} the product of the frequencies of the decode letters which we get if the key is~$\beta$.

%
%
Since the \textit{a priori} probabilities of the keys are all equal we may say that the \textit{a posteriori} probabilities are in the ratio
$\prod_{i} \mathcal{P}_{\alpha_{i} - \beta + 1}$ \textit{i.e.} in the ratio $\prod_{i} 26 \mathcal{P}_{\alpha_{i} - \beta + 1} $ which is more convenient for calculation. The final value for the probability is then
\large{
\[
	\frac{\prod\limits_{i} 26 \mathcal{P}_{\alpha_{i} - \beta + 1} }{\sum\limits_{\beta} \prod\limits_{i} 26 \mathcal{P}_{\alpha_{i} - \beta + 1 } }.
\]}
\normalsize
The calculation of the product $\prod_{i} 26 \mathcal{P}_{\alpha_{i} - \beta + 1 }$
may be done by the method recommended before for  
\[
	\prod_{i} \frac{25 \mathcal{P}_{\alpha_{i} - \beta + 1 }}{ 1 - \mathcal{P}_{ \alpha_{i} - \beta + 1 }} .
\]
$\bigl($The table in Fig 3 was in fact made up for $\prod_{i} 26 \mathcal{P}_{ \alpha_{i} - \beta + 1 } $.
The differences between the two tables would of course be rather slight$\bigr)$. The new result is more accurate than the old because of the independence assumption in the original result.

If we only want to know the ratios  of the probabilities of the various keys there is no need to calculate the denominator $\sum_{\beta}\prod_{i} 26  \mathcal{P}_{ \alpha_{i} - \beta + 1 }$.
This denominator has however another importance: it gives us some evidence about other assumptions, such as that the cipher is Vigen\`{e}re, and that the period is 10. This aspect will be dealt with later (p.  )\footnote{ Forward reference left unresolved in the manuscript.}.	
%
%
\section{A letter subtractor problem}
A substitution with the period $91 \times 95 \times 99$ is obtained by superimposing three substitutions of periods 91, 95, and 99, each substitution being a Vigen\`{e}re composed of slides of 0, 1, 2, 3, 4, 5, 6, 7, 8, or 9.\footnote{ Equivalent to keys A to J.} The three substitutions are known in detail, but we do not know for any given message at what point in the complete substitution to begin. For many messages however we can provide a more or less probable crib. How can we test the probability of a crib before attempting to solve it? It may be assumed that approximately equal numbers of slides 0, 1,~\dots,~ 9 occur in each substitution. 

The principle of the calculation is that owing to the way in which the substitution is built up, not all slides are equally frequent, \textit{e.g.} a slide of 25 can only be the sum of slides of 9, 8, and 8, or 9, 9, and 7 whilst a slide of 15 can be any of the following

\begin{center}
	\begin{tabular}{llll}
		9,6,0		& 8,7,0	&7,7,1 	&6,6,3 \\
		9,5,1 	& 8,6,1	&7,6,2	&6,5,4 \\
		9,4,2		& 8,5,2	&7,5,3 \\
		9,3,3		&8,4,3	&7,4,4
	\end{tabular}
\end{center}
A crib will therefore, other things being equal, be more likely if it requires a slide of 15 than if it requires a slide of 25. The problem is to make the best use of this principle, by determining the probability of the crib with reasonable accuracy, but without spending long over it.

%
%
We have to find the probability of getting a given slide. To do this we can apply several methods.
\begin{enumerate}[(a)]
	\item We can produce a long stretch of key by addition and take a count of the resulting slides. This is obviously a very general method, and requires no special mathematical technique. It may be rather laborious, but by interpreting a small count with common sense one can probably get quite good results.
	\item There are 1000 possible combinations of slides all equally likely \textit{viz.} 000, 001, \dots, 999. We can add up the digits in these and take the remainder on division by 26, and then count the number of combinations giving each of the possible remainders.
	\item We can make use of a trick which might appear to be rather special, but is really applicable to a multitude of problems. Consider the expression
\[
	f(x) = \left(1 + x +x^{2} + \cdots + x^{9} \right)^{3}.
\]
For each possible way of expressing a number \textit{n} as the sum of three numbers 0, \dots, 9, say $n = m_{1} + m_{2} + m_{3}$, there is a term $x^{m_{1}}x^{m_{2}}x^{m_{3}}$ in $f(x)$, $x^{m_{1}}$ coming out of the first factor, $x^{m_{2}}$ out of the second, and $x^{m_{3}}$ out of the third. Hence the number of ways of expressing \textit{n} in the form $n = m_{1} + m_{2} + m_{3}$, is the coefficient of $x^{n}$ in $f(x)$ \textit{i.e.} in
\[
	\frac{\left( 1 - x^{10} \right)^{3}}{\left( 1 - x\right) ^{3}},
\]
or in
\[
	\left( 1 - 3x^{10} + 3x^{20} - x^{30} \right) \left( 1 - x \right) ^{-3}.
\]
%
%
Expanding $\left( 1 - x \right)^{-3}$ by the binomial theorem
\begin{multline*}
	\left( 1 - x \right)^{-3} = 1 + 3x +  6x^{2} + 10x^{3} + 15x^{4} + 21x^{5} + 28x^{6} + 36x^{7} \\
	+ 45x^{8} + 55x^{9} + 66x^{10} + 78x^{11} + 91x^{12} + 105x^{13} \\
	+ 120x^{14} + 136x^{15} + 153x^{16} + 171x^{17} + 190x^{18} \\
	+ 210x^{19} + 231x^{20} + 253x^{21} + 276x^{22} + 300x^{23} \\
	+ 325x^{24} + 351x^{25} + 378x^{26} + 406x^{27} + 435x^{28} + \dots .
\end{multline*}
Now multiply by $1 - 3x^{10} + 3x^{20} - x^{30}$ and we get
\begin{multline*}
	f(x) = 1 + 3x +  6x^{2} + 10x^{3} + 15x^{4} + 21x^{5} + 28x^{6} + 36x^{7} \\
	+ 45x^{8} + 55x^{9} + 63x^{10} + 69x^{11} + 73x^{12} + 75x^{13} \\
	  + 75x^{14} + 73x^{15} + 69x^{16} + 63x^{17} + 55x^{18} \\
	  + 45x^{19} + 36x^{20} + 28x^{21} + 21x^{22} + 15x^{23} \\
	  + 10x^{24} + 6x^{25} + 3x^{26} + x^{27} \qquad \qquad
\end{multline*}
This means to say that the chances of getting totals 0, 1, 2, \dots \, are in the ratio 1, 3, 6, 10, \dots \, The chances of getting remainders of 0, 1, 2, \dots \, on division by~26 are in the ration 4, 4, 6, 10, 15, \dots  \, To get true probabilities these must be divided by their total which is conveniently 1000.
	\item There are two other methods, both connected with the last method but not relying so much on the special features of the problem. They will be discussed later.\footnote{ No such discussion appears in the manuscript.}
\end{enumerate}

Suppose then that the probabilities have been calculated by one method or the other (as in fact we have done under (c)). We can then estimate the values of cribs. Let us suppose that a possible crib for a message beginning \texttt{MVHWUSXOWBVMMK} was \texttt{AMBASSADOR} so that the slides were 12, 9, 6, 22, 2, 0, 23, 11, 14. The slide of 12 gives us some slight evidence in favour of the crib being right for slides of 12 occur with frequency 0.073 with right cribs, whilst with wrong cribs they occur with frequency only 1/26. 
%
%
The factor in favour of the crib is therefore $26 \times 0.073$ or about 1.9. A similar calculation may be made for each of the slides, but of course the work may be greatly speeded up by having the values of the factors 26 $C_{s}/1000$ in half decibans tabulated: here $C_{s}$ is the coefficient of $x^{s}$ in the above polynomial $f(x)$. The table is given below (Fig 6)
\begin{center}
	\begin{tabular}{rrr}
		1	& 0 	&-20  \\
		2	&25	&-16  \\
		3	&24	&-12  \\
		4	&23	&-8  \\
		5	&22	&-6  \\
		6	&21	&-3  \\
		7	&20	&-1  \\
		8	&19	&1  \\
		9	&18	&3  \\
		10	&17	&4  \\
		11	&16	&5  \\
		12	&15	&6  \\
		13	&14	&6
	\end{tabular}
\end{center}
\begin{center}
	\textsc{Figure 6.} Scores in half decibans of the various slides.
\end{center}

Evaluating this crib by means of this table we score
\[
	6 + 3 - 3 - 6 -16 - 20 - 8 + 5 + 6 \ ( \ = -33\ ),
\]
\textit{i.e.} the crib is worse by a factor of $10^{-33/20}$ than it was before \textit{e.g.} if the \textit{a priori} odds of the crib were 2:1 against it becomes 98:1 against. This crib was in fact made up at random \textit{i.e.} the letters of the cipher text were chosen at random. 

%
%
Now let us take one made up correctly, \textit{i.e.} really enciphered by the method in question, but with a random chosen key.

\Large
\begin{bfseries}
\begin{center}
	\begin{ttfamily}
		\begin{tabular}{c c c c c c c c c c}
			N &Y &X &L &N &X &I &Q &H &H  \\
			A &M &B &A &S &S &A &D &O &R  \\
			13 &12 &22 &11 &21 &5 &8 &13 &19 &16		
		\end{tabular}
	\end{ttfamily}
\end{center}
\end{bfseries}
\normalsize
\begin{center}
	\begin{itshape}
		(slides)	
	\end{itshape}
\end{center}
This scores 15 so that if it were originally 2:1 against, it now becomes nearly 3:1~on.

Having decided on a crib the natural way to test it is to have a catalogue of the positions in which a given series of slides is obtained if the 91 period component is omitted. We make 91 different hypotheses as to this third component, draw an inference as to what is the part of the slide arising from the components of periods 95 and 99 combined. This we look up in the catalogue. This process is fairly lengthy, and as the scoring of the crib takes only a minute it is certainly worth doing.
%
\section{Theory of repeats}
Suppose we have a cipher in which there are several very long series of substitutions which can be used for enciphering a message, but that one  may sometimes get two messages enciphered with the same series of substitutions (or possibly, the series of substitutions for one message being those for another with some at the beginning omitted). In such a case let us say that the messages \textit{fit}, or that they fit at such and such a distance, the distance being the number of substitutions which have to be omitted from the one series to obtain the other series. One will frequently want to know whether two messages fit or not, and we may find some evidence about this by examining the repeats between them. 

By the repeats between them I mean this. One writes out the cipher texts of the two messages with the letters which are thought to have been enciphered with the same substitution under one another. One then writes under these messages a series of letters \texttt{O} and \texttt{X}, an \texttt{O} being written where the cipher texts differ and an \texttt{X} where they agree. The series of letters \texttt{O} and \texttt{X} will begin where the second message begins and end where the first to end ends. This series of letters \texttt{O} and \texttt{X} may be called the repetition figure. It may be completed by adding at the ends an indication of how many letters there are which do not overlap, and which message they belong to. 

As an example:
\large
\texttt{
	\begin{center}
		GFRLIKQGVBMILAFIXMMOROGBYSKYXDAZCHMUMRKBZLDLDDOHCMVTIPRSD\\
		VLOVDYQCEJSOPYGBMBKYXDAZNBFIOPTFCXDOD\; \; \;    \\
		$^{8}$XOOOOOOOOOOXOOXXOOXXXXXXOOOOOOOOOOXOX$^{11}$\, \, \,
	\end{center}}
\normalsize

%
%
On the whole one expects that a fit is more likely to be right the more letters \texttt{X} there are in the repetition figure, and that long series of letters \texttt{X} are especially desirable. This is because it would not be very unusual for two fairly common words to lie directly under one another when the clear texts are written out, thus
\large
\texttt{
	\begin{center}
		THEMAINCONVOYWILLARRIVE \dots \\
		\qquad \ ALLCONVOYSMUSTREPORT \dots \\
		\qquad  XOOXXXXXXOOOOOXOOOO \dots
	\end{center}}
\normalsize
If the corresponding cipher texts really fit, \textit{i.e.} if the letters in the same column are enciphered with the same substitution, then the condition for an \texttt{X} in the repetition figure of the cipher texts is that there be an \texttt{X} in the repetition figure of the corresponding clear text. Now series of several consecutive letters \texttt{X} can occur quite easily as above by two identical words coming under one another, or by such combinations as
\large
\texttt{
	\begin{center}
		ITISEASIERTOTEACHTHANALGEBRA \dots \\
		\quad \: THERAINWASSUCHTHATHECOULD \dots \\
		\quad \: OOOOOOOOOOOOXXXXXOOOOOOOO \dots
	\end{center}}
\normalsize
if the messages really fit, but if not they can only occur by complete coincidence. One therefore tends to believe that there is a fit when one gets such series of letters~\texttt{X}. As regards single cases of \texttt{X} the value of them is not so clear, but one can see that if $\mathcal{P}_{\alpha}$ is the frequency of letters $\alpha$ in plain language then the frequency of letters \texttt{X} as a whole in comparison of plain language with plain language is $\sum_{\alpha}\mathcal{P}_{\alpha}^{\;2}$, whilst for wrong fits of cipher text it is 1/26 which is necessarily less. Given a sufficiently long repetition figure one should therefore be able to tell whether it is a fit or not simply by counting the letters \texttt{X} and \texttt{O}.

%
%
So much is well known. The real point of this section is to show these ideas can be developed into an accurate method of estimating the probabilities of fits.

\subsection{Simple form of theory}
The complete theory takes account of the various possible lengths of repeat. As this theory is somewhat complicated it will be as well to give first two simplified forms of the theory. In both cases the simplification arises by neglecting a part of the evidence. In the first simplified form of theory we neglect all evidence except the number of letters \texttt{X} and the number of letters \texttt{O}. In the other simplified form the evidence is the number of series of (say) four consecutive letters \texttt{X} in a repetition figure. 

When our evidence is just the number of times \texttt{X} occurs in the repetition figure, (\textit{n} let us say) and the length of the repetition figure (\textit{N} say), then the factor in favour of the fit is 
\[
	\frac{\textit{Probability of a right repetition figure of length N and n occurrences of \texttt{X}}}{\textit{Probability of a wrong repetition figure of length N having n occurrences of \texttt{X}}}.
\]
As an approximation we may assume that the numerator of this expression has the same value as if the right repetition figures were produced letter by letter by independent random choices, with a certain fixed probability of getting an \texttt{X} at each stage. This probability will have to be $\beta = \sum_{\alpha}\mathcal{P} _{\alpha}^{\;2}$. The numerator is then
%
%
\small{
\begin{multline*}
	\left( \textit{Number of repetition patterns with length N and n occurrences of \texttt{X}} \right) \\
	\times \left( \textit{Probability of getting a given such repetition pattern by the process just mentioned} \right),
\end{multline*}}
\normalsize
which we may write as $R(N,n)Q(N,n)$. Now let us denote by $y_{i}$ the $i$th symbol of the given repetition pattern and put $\tau_{x} = \beta$ and $\tau_{0} = 1 - \beta$. Then $Q(N,n)$, the probability of getting the repetition pattern is $\prod_{i=1}^{N}\tau_{y_{i}}$ which simplifies to $\beta^{n}(1 - \beta)^{N-n}$. We may do a similar calculation for the denominator, but here we must take $\beta = 1/26$ since all letters occur equally frequently in the cipher. The denominator is then
\[
	R(N,n) \left( \frac{1}{26} \right)^{n} \left( \frac{25}{26} \right) ^{N-n}.
\]
In dividing to find the factor for the fit $R(N, n)$ cancels out, leaving
\[
	\left( 26\beta \right) ^{n} \left( \frac{26}{25} \left( 1 - \beta \right) \right)^{N - n}.
\]
In other words we score a factor of $26 \beta$ for an \texttt{X} and a factor of $(26/25)(1 - \beta)$ for an \texttt{O}. More convenient is to regard it as $10 \log_{10} \bigl[(25 \beta) /(1 - \beta)\bigr]$ decibans for an \texttt{X} and  $10 \log_{10}[(26/25) /(1 - \beta)]$ per unit length of repetition figure (\textit{per unit overlap}).

An alternative argument, leading to the same result, runs as follows. Having decided to neglect all evidence except the overlap and the number of repeats we pretend that nothing else matters, \textit{i.e.} that the form of the figure is irrelevant. In this case we can regard each letter of the repetition figure as independent evidence about the fit. If we get an \texttt{X} the factor for the fit is
\[
	\frac{\textit{Probability of getting an \texttt{X} if the fit is right}}{\textit{Probability of getting an \texttt{X} if the fit is wrong}},
\]
\textit{i.e.} $\beta/(1/26)$.
Similarly the factor for an \texttt{O} is $(1 - \beta)/(25/26)$.

%
%
In either form of argument it is unnecessary to calculate the number  $R(N,n)$. In this particular case there is no particular difficulty about about it: it is the binomial coefficient. In some similar problems this cancelling out is a great boon, as we might not be able to find any simple form for the factor which cancels. The cancelling out is a normal feature of this kind of problem, and it seems quite natural that it should happen when we think of the second form of argument in which we think of the evidence as consisting of a number of independent parts.

The device of assuming, as we have done here, that the evidence which is not available is irrelevant can often be used and usually leads to good results. It is of course not supposed that the evidence really is irrelevant, but only that the error resulting from the assumption when used in this kind of way is likely to be small.

\subsection{Second simplified form of theory}
In the second simplified form of theory we take as our evidence that a particular part of the repetition figure is \texttt{OXXXXO} (say, or alternatively \texttt{OXXXXXO} say). The factor is then
\[
	\frac{\textit{Frequency of \texttt{OXXXXO} in right repetition figures}}{\textit{Frequency of \texttt{OXXXXO} in wrong repetition figures}}.
\]
The denominator is
\[
	\left( \frac{1}{26} \right) ^{4} \left( \frac{25}{26} \right)^{2},
\]
and the numerator may be estimated by taking a sample of language hexagrams and counting the number of pairs that have the repetition figure \texttt{OXXXXO}. The expectation of the number of such pairs is the sum for all pairs of the probabilities of those pairs having the desired repetition figure \textit{i.e.} is the number of such pairs (\textit{viz} $N(N-1)/2$ where $N$ is the size of the sample) multiplied by the frequency of \texttt{OXXXXO} repetition figures. This frequency may therefore be obtained by division if we equate the expected number of these repetition figures to the actual number.
%
%
\newpage
\subsection{General form of theory}
It is not of course possible to have statistics of every conceivable repetition figure. We must make some assumptions to reduce the variety that need to be considered. The following assumption is theoretically very convenient, and also appears to be a very good approximation.

\textit{The probability of repeats at two points known to be separated by a point where there is known to be no repeat are independent}.

We may also assume that the probability of a repeat is independent of anything but the repetition figure in this neighbourhood. (We may however as a refinement produce different positions in a message). We can therefore think of repetition figures as being produced by selecting the symbols of the figure consecutively, the probability of getting an \texttt{X} at each stage being determined by the repetition figure from the point in question back as far as the last \texttt{O}. Sometimes this will take us back as far as the beginning of the message, and will include the number telling us how many more letters there are which do not repeat at all. We need in practice only distinguish two cases, where this number is  0 and when it is more. We may also neglect the question as to which message occurs first. We therefore have to distinguish the following cases
%
%
\begin{center}
	\begin{tabular}{lrlrlr}
		\texttt{O} &$a_{0}$ &some &$b_{0}$ &none &$c_{0}$  \\
		\texttt{OX} &$a_{1}$ &some \texttt{X} &$b_{1}$ &none \texttt{X} &$c_{1}$  \\
		\texttt{OXX} &$a_{2}$ &some \texttt{XX} &$b_{2}$ &none \texttt{XX} &$c_{2}$ \\
		\texttt{OXXX} &$a_{3}$ &some  \texttt{XXX} &$b_{3}$ &none \texttt{XXX} &$c_{3}$\\
		\dots &		&\dots	  &		 &\dots		
	\end{tabular}
\end{center}
\normalsize
The entries $a_{0}, a_{1}, b_{0}, etc.$ opposite the repetition figures are the notations we are adopting for the probability of getting another \texttt{X} following such a figure. Strictly speaking we should also bring in a notation for the probability of the message coming to an end after any given repetition figure. As the repeats at the end of a comparison do not appear to behave very differently from those in the main part of the message I shall neglect this complication by assuming that the probability of getting an \texttt{O} added to the probability of getting an \texttt{X} is 1, and that afterwards one cuts off the end of the series arbitrarily.

Let us calculate the factor for the repeat figure\footnote{ In the manuscript, Turing squeezes the figure into three lines by spilling into the margins and use of pen and ink. The typeset equivalent is unreadable, so the figure has been split into a left and right components.

Reassemble as: \texttt{none X X X X O | O | O | X O  | X X X O | O | X X | some}}
%
%
\begin{figure}[ht]
\raggedright\includegraphics[scale=0.6]{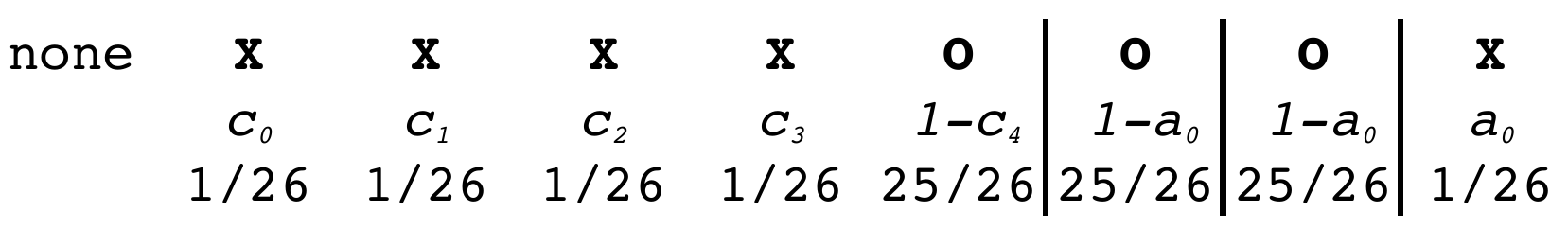}
\end{figure}
\begin{figure}[ht]
\raggedleft\includegraphics[scale=0.6]{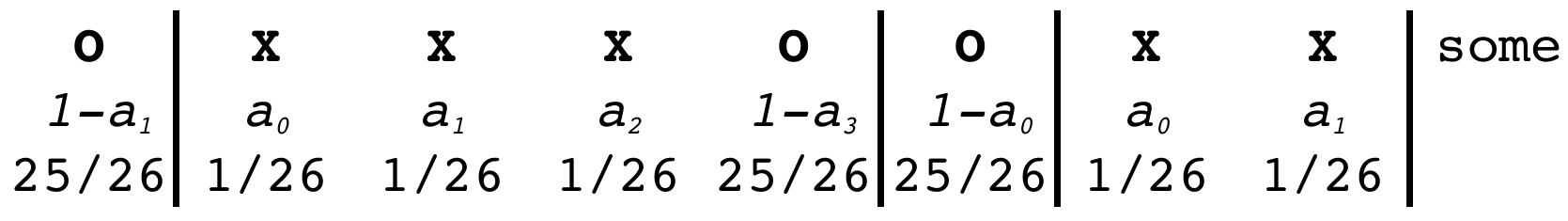}
\end{figure}

Underneath each symbol has been written the probability that one would get that symbol, knowing the ones which precede, both for the case of a right and of a wrong repetition figure. The factor for the fit is the product of the first row divided by the product of the second. It is convenient to split this up as indicated by the vertical lines into the product of
%
%
\begin{align*}
	\quad \quad &\frac{c_{0}c_{1}c_{2}c_{3} \left( 1 - c_{4} \right)}{(1/26)^{4} \times \left(25/26\right)}, \\
	\quad \quad & \frac{1 - a_{0}}{\left(25/26\right)}, \qquad \text{- occurring three times}, \\
	\quad \quad & \frac{a_{0} (1 - a_{1} )}{\left(1/26\right) \times \left(25/26\right)}, \\
	\quad \quad & \frac{a_{0}a_{1}a_{2} \left( 1 - a_{3} \right)}{(1/26)^{3} \times \left(25/26\right)}, \\
	\quad \quad & \frac{a_{0} a_{1}}{(1/26)^{2}},
\end{align*}

and this product may be put into the form of the product of 
\[
	\frac{c_{0}c_{1}c_{2}c_{3} \left( 1 - c_{4} \right)}{\left( 1/26 \right)^{4} \times \left( 25/26 \right)} \times  \left( \frac{1 - a_{0}}{\left(25/26\right)} \right)^{-5},
\]
\begin{flushright}
- which we call the factor for an \\ initial tetragramme repeat level,
\[
	\frac{a_{0} (1 - a_{1} )}{\left( 1/26 \right) \times \left( 25/26 \right)} \times   \left( \frac{1 - a_{0}}{\left(25/26\right)} \right) ^{-2},
\]
- the factor for a single repeat,
\[
	\frac{a_{0}a_{1}a_{2} \left( 1 - a_{3} \right)}{\left( 1/26 \right)^{3} \times \left( 25/26 \right)}  \times  \left( \frac{1 - a_{0}}{\left(25/26\right)} \right) ^{-4},
\]
- the factor for a trigramme,
\[
	\frac{1 - a_{0}}{1 - a_{2}},
\]
- the correction for a final bigramme,
\[
	\left( \frac{1 - a_{0}}{\left(25/26\right)} \right) ^{16},
\]
- the factor for an overlap of 16,
\[
	\frac{a_{0}a_{1} \left( 1 - a_{2} \right)}{ \left( 1/26  \right)^{2} \times \left( 25/26 \right)} \times  \left( \frac{1 - a_{0}}{\left(25/26\right)} \right) ^{-3},
\]
- the factor for a trigramme.
\end{flushright}
	
%
%
We shall neglect the correction for a final bigramme (or whatever it may be). It is in any case rather small, and vanishes if the repetition figure ends with \texttt{O}; also with our conventions the whole question of the ends of repetition figures has been left rather in doubt.

Now let us put\footnote{ The manuscript has a pencilled note beside $k_{r}$ indicating it is to be read as $k_{r+1}$. We presume that this also means that $j_{r}$ should be $j_{r+1}$, and $i_{r}$ should be $i_{r+1}$, However, these are not indicated and no changes are made in the subsequent text. We leave the text unchanged}
\begin{gather*}
	a_{0} a_{1} \dots a_{r} ( 1 - a_{r+1} ) = k_{r},  \\
	b_{0} b_{1} \dots b_{r} ( 1 - b_{r+1} ) = j_{r},  \\
	c_{0} c_{1} \dots c_{r} ( 1 - c_{r+1} ) = i_{r}.
\end{gather*} 
The values of the $i_{r}$ can be obtained as follows. We take a number of plain language messages and leave out two or three words at the beginning. Then combine the messages to form one long message; this message may be made to \textit{eat its own tail} \textit{i.e.} it may be written round a circle. If the message were compared with itself in every possible position, except level, we should expect to get repetition figures which when divided up as shown by vertical lines after each \texttt{O}, containing $(N(N-1)/2) k_{r} \  (= N_{r})$ parts which consist of \textit{r} symbols \texttt{O}, or as we may say $N_{r}$ \textit{actual r-gramme repeats}, where \textit{h} is the probability of an \texttt{O} . 

The values of $N_{r}$ can be calculated given the \textit{apparent number of \textit{r}-gramme repeats} $M_{r}$ for each $r$. This apparent number of \textit{r}-gramme repeats is the number of series of \textit{r} consecutive symbols \texttt{X} in the repetition figures regardless of what precedes or follows the series.

By considering the ways in which an actual repeat can give rise to the apparent repeat of various lengths we see that
\[
	M_{r} = N_{r} + 2N_{r+1} + 3N_{r+2} + \dotsc ,
\]
and therefore
\[
	M_{r} - M_{r+1} = N_{r} + N_{r+1} + N_{r+2} + \dotsc ,
\]
and
\[
	\left( M_{r} - M_{r+1} \right) - \left( M_{r+1} - M_{r+2} \right) = N_{r}.
\]
%
%
The calculation of $j_{r}$ may perhaps best be done by comparing the beginners of a number of messages with the long circular message, and the values of $i_{r}$ by comparing the beginners among themselves. A similar technique of actual and apparent numbers of repeats can be used. I shall not go into this in detail. The formulae required may now be assembled.
\begin{align*}
	\mu_{r} &= \text{decibanage for an \textit{r}-gramme repeat}, \\
	\gamma &= \text{negative decibanage for unit overlap}, \\
	S_{\beta, r} &= \text{number of occurrences in the statistics of the  \textit{r}-gramme $\beta$}, \\
	N &= \text{total number of letters in the statistics}. 
\end{align*}
Then if
\begin{align*}
	M_{r} &= \sum_{\beta}  \frac{S_{\beta, r} \left( S_{\beta, r} - 1 \right)}{2},  \\
	N_{r} &= M_{r} - 2M_{r+1} + M_{r+2}, \\
	L &= \frac{N(N -1)}{2},  \\
	k_{r} &= \frac{N_{r}}{Lh}.
\end{align*}
$h$ may be calculated as follows. From the identity
\[
	\left( 1 - a_{0} \right) + a_{0} \left( 1 - a_{1} \right) + a_{0}a_{1} \left( 1 - a_{2} \right) + \dots = 1, 
\]
we get
\[
	k_{0} + k_{1} + k_{2} + \dots = 1,
\]
\[
	\therefore \qquad \frac{L - M_{1}}{Lh} = 1,
\]
\[
	\left( 1 - a_{0} \right) = k_{0} = \frac{N_{0}}{L - M_{1}} = \frac{L- 2M_{1} + M_{2}}{L - M_{1}},
\]
\[
	\mu_{r} = 10 \log_{10} \left( \frac{26^{r+1} k_{r}}{25} \right) + \left( r + 1 \right)\nu,
\]
\[
	\nu = - 10 \log_{10} \left( \frac{26 ( 1 - a_{0})}{25} \right).
\]

%
%
\section{Transposition ciphers}
\subsection{A probability problem}
In making calculations about substitution ciphers we have often found it useful to treat the plain language as if it were produced by independent choices for the letters, using certain fixed frequencies with which the letters are chosen. Our method for Vigen\`{e}re and one of the simplified forms of repeat theory could be based on this sort of assumption.  With a transposition cipher however such an assumption would be useless or worse than useless, for it would result in the conclusion that all transpositions were equally likely. We have therefore to take a slightly less crude assumption, and the one which suggests itself is that the letters forming the plain language are chosen consecutively, the probability of getting a particular letter depending only on what the letter is and what the preceding letter was. It is easily verified the if $\mathcal{P}_{\alpha \beta}$ is the proportion of bigrammes $\alpha \beta$ in plain language and $\mathcal{P}_{\alpha}$ the frequency of the letter $\alpha$ then the probability $q_{\alpha \beta}$ of a letter $\beta$ following an $\alpha$ is $\mathcal{P}_{\alpha \beta}/\mathcal{P_{\alpha}}$. The probability of a piece of plain language of length $ L $ letters saying $\alpha_{1} \alpha_{2} \dots \alpha_{L}$ is then
\[
	\mathcal{P}_{\alpha_{1}} \times  q_{\alpha_{1}\alpha_{2}} \times q_{\alpha_{2}\alpha_{3}} \times q_{\alpha_{3}\alpha_{4}}   \times \dots \times q_{\alpha_{(L-1)}\alpha_{L}},
\]
which may also be written as
\[
	\mathcal{J} \left( \alpha_{1}, \dotsc, \alpha_{L} \right).
\]
We may also calculate the probability of a given piece of plain language having certain given letters in given places, the remainder of the message being unspecified. The probability is given by
%
%
\[
	\sum \left( \xi_{1} , \dotsc , \xi_{L} \text{ consistent with data } \right) \mathcal{J}  \left( \xi_{1},  \dotsc , \xi_{L} \right),
\]
and if the data is that the known letters are 
\begin{equation}\label{E:letterD}
	\underset{ n_{1} \text{ dots }}{\cdots} \beta_{1} \quad \underset{ n_{2} \text{ dots } }{\cdots} \beta_{2} \quad \cdots \quad\cdots \beta_{r-1} \quad \underset{ n_{r} \text{ dots }}{\cdots} \beta_{r} \quad \cdots,
\end{equation}
it is approximately\footnote{ The manuscript has as the first term $\prod_{r} \beta_{r}$, a pencilled annotation indicates that  the $\beta_{r}$ is to be read as $\mathcal{P}_{\beta_{r}}$. This substitution has been made in the text.}
\begin{equation}\label{E:letterA}
	\prod_{r} \mathcal{P}_{\beta_{r}} \cdot \prod_{n_{r+1} = 0} \frac{\mathcal{P}_{\beta_{r}\beta_{r+1}}}{\mathcal{P}_{\beta{r}}\mathcal{P}_{\beta_{r+1}}}.
\end{equation}
 A more or less rigorous deduction of this approximation from the assumptions above is given at the end of the section. For the present let us see how it can be applied. If we have two theories about the transposition of which the one requires the above pattern of letters, and the other brings the same letters in to positions in which no two of them are consecutive, then the factor in favour of the first as compared with the second is
 \[
 	\prod_{n_{r+1} = 0} \frac{\mathcal{P}_{\beta_{r}\beta_{r+1}}}{\mathcal{P}_{\beta{r}}\mathcal{P}_{\beta_{r+1}}}.
 \] 
 We can apply this straightforwardly to the case of a simple transposition by columns. 
%
%
The following text is known to be a simple transposition of a certain type of German text with a key length of not more than 15.\footnote{ As for the Vigen\`{e}re problem above, Turing's statement of the ciphertext is slightly different from that which he scores for decryption. The second line in the ciphertext below begins \texttt{NLTS}, however, this changes to \texttt{NITS} in the scoring example in Figure 7 below. See also the notes accompanying the cleartext.}
\vspace{0.1in}
\large
\texttt{
	\begin{flushleft}
		\quad  S A T P T W S F A S T A U T E E A I E U F H W T J T D D G C  \\
		\quad  N L T S E F C U I E B O E Y Q H G T J T E E F I E O R T A R  \\
		\quad  U R N L N N N N A I E O T U S H L E S B F B R N D X G N J H  \\
		\quad  U A N W R	
	\end{flushleft}}
\normalsize
\vspace{0.1in}

To solve this transposition, we may try comparing the first six letters of \linebreak \texttt{S~A~T~P~T~W} which we know form part of one column with each other series of six letters in the message, for we know that one such comparison will give entirely bigrammes occurring in the decode. We may try first
\vspace{0.1in}
\large
\texttt{
	\begin{center}
		S F  \\
		A A  \\
		T S  \\
		P T  \\ 
		T A  \\
		W U
	\end{center}}
\normalsize
\vspace{0.1in}

The factor for a transposition which brings these letters together, as compared with one which leaves them apart is
\[
	\frac{\mathcal{P}_{SF}}{\mathcal{P}_{S}\mathcal{P}_{F}} \times \frac{\mathcal{P}_{AA}}{\mathcal{P}_{A}\mathcal{P}_{A}} \times \dots \times \frac{\mathcal{P}_{WU}}{\mathcal{P}_{W}\mathcal{P}_{U}}.
\]
By using a table of values of 
\[
	20 \log_{10} \left( \frac{\mathcal{P}_{\alpha \beta}}{\mathcal{P}_{\alpha}\mathcal{P}_{\beta}} \right),
\]
made up for the type of traffic in question, and given to the nearest integer (table of values of $\mathcal{P}_{\alpha \beta}/(\mathcal{P}_{\alpha}\mathcal{P}_{\beta})$ expressed in half-decibans) we get the product by addition. Such a table is shown in Fig 6. The scores for this particular columns are \texttt{ SF -7, AA -7, TS -2, PT -10, TA -3, WU -13}, totalling -36. If we consider this combination as \textit{a priori} about 100:1 against (there are 95 letters in the message) it is \textit{a posteriori} about 3000:1 against.

%
%
\newpage
\begin{landscape}
\begin{figure}[ht]
\centering\includegraphics[scale=0.85]{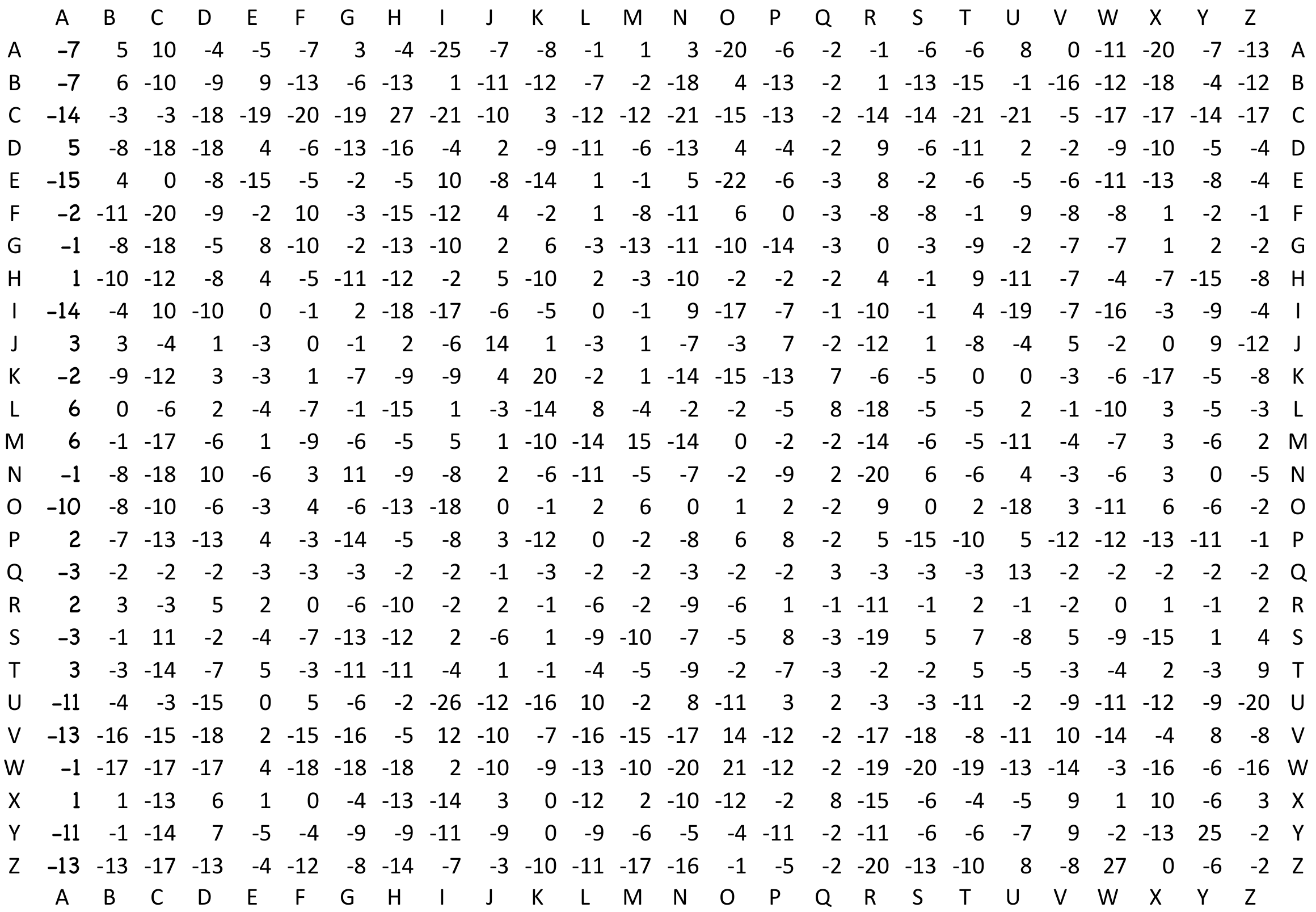}
\caption{Exclusive bigramme scores in half decibans, \textit{i.e.} $20 \log_{10} \left( \frac{\mathcal{P}_{\alpha \beta}}{\mathcal{P}_{\alpha}\mathcal{P}_{\beta}} \right)$, for a certain kind of German traffic.}\label{Fi:FIG6}
\end{figure}
\end{landscape}
%
%
%
%
\newpage
\begin{figure}[ht]
\centering\includegraphics[scale=0.79]{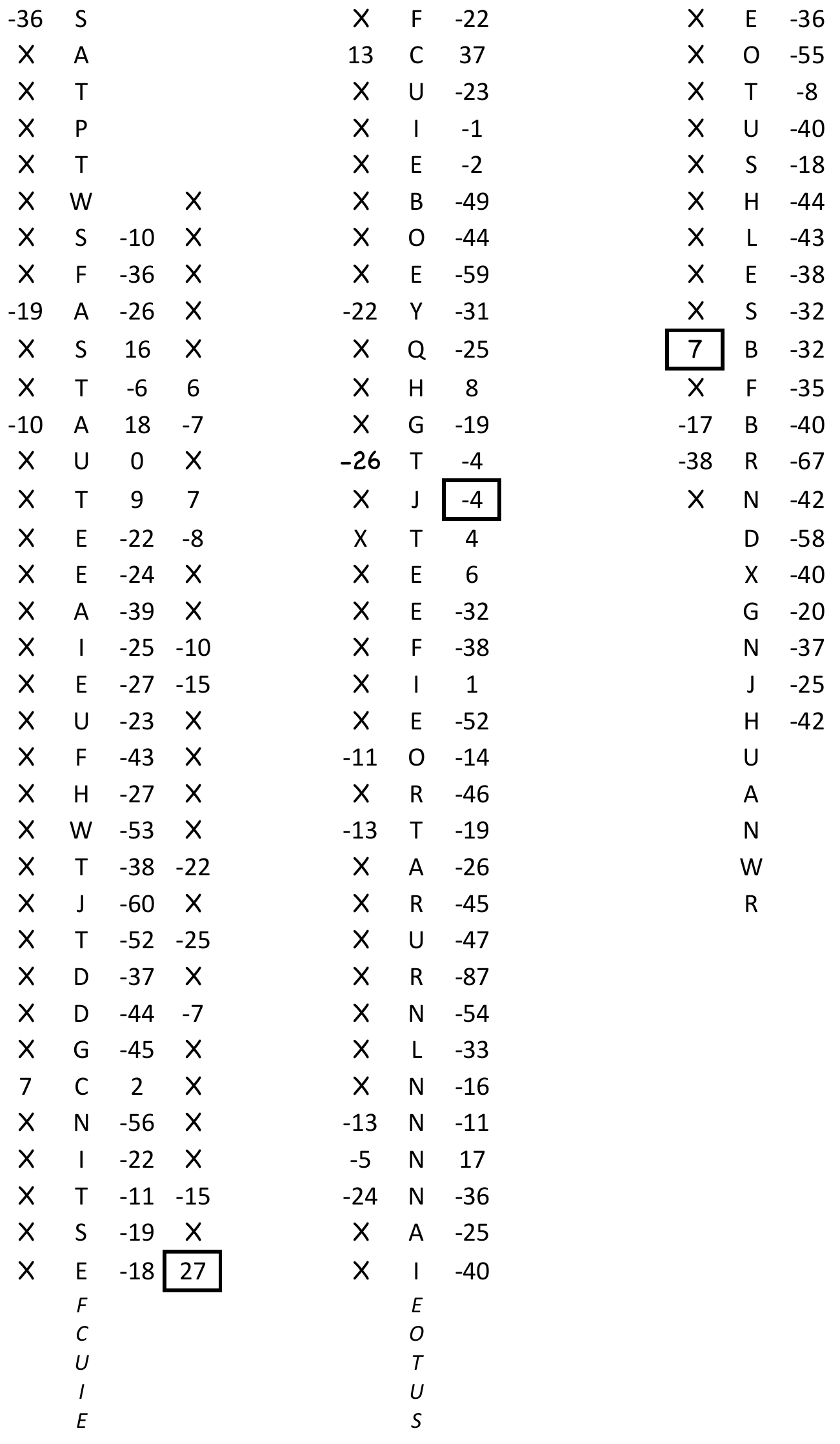}
\caption{Scoring the matching of columns in a simple transposition. Correct matchings noted.}\label{Fi:FIG7}
\end{figure}

\newpage
%
%
Similar scoring may be done for every possible comparison of \texttt{S A T P~T~W} with six consecutive letters of the message. The comparison may be made both with \texttt{S A T P T W} as earlier and as later column; one may also use the last six letters of the message \texttt{H U A N W R}. 

The results of doing this are shown in Fig 7. The message has been written out vertically. The first columns of figures after the message gives the score for \linebreak\texttt{S A T P T W} as earlier column, entered against the first letter of the later column, \textit{e.g.}  the -36 as calculated above gets entered against the  \texttt{F} of \texttt{F A S T A U}. The second column after the message consists of the scores for \texttt{H U A N W R} as first column [and the column before the message gives scores for \texttt{H U A N W R} as second column].\footnote{ \emph{\dots and the column} to end of sentence, has the note in pencil: \emph{I doubt it - S.W.}} One of these columns has been worked out in detail but in the other two crosses have been put in where the scores are very bad.
 
The scores which eventually turned to to be right are ringed. The fourth comparison, which did not have to be done scored very badly \textit{viz.} -27. Amongst the good scores which were wrong there was one score of 37. It was not difficult to see that this one was wrong as most of the score came from \texttt{W O} with requires \texttt{Z} to precede it, and there was no \texttt{Z} in the message. Apart from this fact the comparison was about evens, although if we take into account the fact that there was no better score it would be better.\footnote{ Using Turing's scoring recommendations and a key length of 12 with sequence 5, 11, 8, 7, 3, 10, 6, 12, 9, 4, 1, 2,  a cleartext emerges: \texttt{ BNTO SJJ ALBA RFJ STATT IN OST B HEUTE DEN ETA RUFS PEDUNYAR NACHT FGFNQUUDNUL WICH AHTR X WIESEN WI GEN GRESFOITE TE}. With Turing's original statement of the ciphertext, as noted above, \texttt{GRESFOITE} becomes \texttt{GRESFOLTE}. Turing scores for the \texttt{I} and not for  \texttt{L}, although it makes no differences in the decision to align bigrams.}  [We have already had a case of this kind of thing in connection with Vigen\`{e}re; if the various positions are \textit{a priori} equally likely and the factors are $f_{1}, f_{2}, \dots, f_{N}$ then the value $f_{r}/\sum f_{i}$ for the probability of the $r$th alternative is better than $\left( f_{r}/N \right) / \left( 1 + f_{r}/N \right)$].
%
%
\subsection{The Probability Formula} 
(Semi) Rigorous deduction of the formula (\ref{E:letterA}) on page \pageref{E:letterA}.
(This is  something of a digression).

The probability of a piece of plain language coinciding where necessary with the data (\ref{E:letterD}) on page \pageref{E:letterD} is
\[
	\mathcal{P}_{\beta_{1}}   \mathcal{T}_{n_{2},\beta_{1} \beta_{2}} \mathcal{T}_{n_{3},\beta_{2} \beta_{3}} \dots  \mathcal{T}_{n_{m},\beta_{m-1} \beta_{m}},
\]
where
\[
	 \mathcal{T}_{n,\alpha \beta} \quad  \text{is} \quad \sum_{\eta_{1} \eta_{1} \dots \eta_{n}} q_{\alpha \eta_{1}} q_{\eta_{1} \eta_{2}} \dots  q_{\eta_{n} \beta},
\]
since
\[
	\sum_{\eta_{1} \dots \eta_{n_{1}}} \mathcal{P}_{\eta_{1}} q_{\eta_{1} \eta_{2}} \dots q_{\eta_{n_{1}} \beta_{1}} = \mathcal{P}_{\beta_{1}}.
\]
We can put
\[
	\mathcal{T}_{n,\alpha \beta} = \left( \mathcal{Q}^{n+1} \right)_{\alpha \beta},
\]
where $\mathcal{Q}$ is the matrix whose $\alpha \beta$ coefficient is $q_{\alpha \beta}$.
The formula (\ref{E:letterA}) on page \pageref{E:letterA} would then be accurate if we could say that for $n > 0,$
\[ 
	\left( \mathcal{Q}^{n+1} \right)_{\alpha \beta} = \mathcal{P}_{\beta}.
\]
This is not true, but it is true that except for very special values for $q_{\alpha \beta}$,
\[
	\left( \mathcal{Q}^{n} \right)_{\alpha \beta} \to \mathcal{P}_{\beta}, \text{ as } n \to \infty,
\]
and the convergence is rather rapid. 

To prove this I shall assume that the eigenvalues of $\mathcal{Q}$ are all different in modulus. In this case we can find a matrix $\mathcal{U}$ with unit determinant, such that $\mathcal{U}^{-1}\mathcal{Q \ U}$ is in the diagonal form
\[
	\mathcal{M} = \mathcal{U}^{-1}\mathcal{Q \ U} = 
	\left(
		\begin{array}{ccccc}
		\mu_{1}	& 0			& 0		& \cdots		& 0			\\
		0		& \mu_{2}		& 0		& 			& \vdots		\\
		0		& \ddots		& \ddots	& \ddots		& 0			\\
		\vdots	&			& 0		& \mu_{25}	& 0			\\
		0		& \cdots		& 0		& 0			& \mu_{26}
		\end{array}
	\right),
\]
%
%
since $\mathcal{QU = UM}$ we have
\[
	\sum_{\gamma} q_{\alpha \gamma} u_{\gamma \beta} = \sum_{\xi} u_{\alpha \xi} m_{\xi \beta},
\]
\textit{i.e.}
\[
	\sum_{\gamma} q_{\alpha \gamma} u_{\gamma \beta} = \mu_{\beta} \mu_{\alpha \beta}.
\]
That is, for each $\beta, u_{\alpha \beta}$ provides a solution of
\begin{equation}\label{E:letterE}
	\sum_{\gamma} q_{\alpha \gamma} l_{\gamma} = \mu l_{\alpha},
\end{equation}
with $\mu = \mu_{\beta}$. Conversely if we have any solution of (\ref{E:letterE}) then $\mu= \mu_{\vartheta} , l_{\alpha} = ku_{\alpha \vartheta	}$ for some $k, \vartheta$ and all $\alpha$, for as $\mathcal{U}$ is non singular we can find numbers $c_{\gamma}$ such that
\[
	l_{\alpha} = \sum_{\gamma} u_{\alpha \gamma} c_{\gamma} \text{ for all } \alpha,
\]
and then substituting in (\ref{E:letterE}) we get
\[
	\sum_{\gamma, \delta} q_{\alpha \gamma} u_{\gamma \delta} c_{\delta} = \mu \sum_{\delta} u_{\alpha \delta} c_{\delta},
\]
\textit{i.e.}
\[
	\sum_{\delta} \left( \mu_{\delta} c_{\delta} - \mu c_{\delta} \right) u_{\alpha \delta} = 0.
\]
Which, since $\mathcal{U}$ is nonsingular implies $\mu  = \mu_\delta  \text{ or } c_{\delta} = 0  \text{ for all } \delta$.

As the series $\mu_{1}, \dotsc, \mu_{26}$ are all different there is only one value $\vartheta$ of $\delta$ for which $\mu = \mu_{\delta}$ and so $l_{\alpha} =  c_{\vartheta}u_{\alpha \vartheta}$ for all $\alpha$.
%
%
Now putting $l_{\alpha} = 1$ for all $\alpha$ we see that one member of the series $\mu_{1}, \dotsc, \mu_{26}$ is 1, for (\ref{E:letterE}) is certainly satisfied. 

I shall prove that the remaining eigenvalues satisfy $\vert \mu \vert \leq 1$. We first prove that if $\mu \neq 1$ then $\sum p_{\alpha} l_{\alpha} = 0$. This follows by multiplying (\ref{E:letterE}) on each side by $\mathcal{P}_{\alpha}$ and summing. Since
\[
	q_{\alpha \beta} = \frac { \mathcal{P}_{\alpha \beta}}{ \mathcal{P}_{\alpha}} \text{ and } \sum_{\alpha}  \mathcal{P}_{\alpha \beta} = \mathcal{P}_{\beta},
\]
we get
\[
	\sum_{\alpha \gamma} q_{\alpha \gamma} l_{\gamma} = \sum p_{\gamma} l_{\gamma} = \mu \sum p_{\alpha} l_{\alpha},
\]
which implies
\[
	\mu = 1 \text{ or } \sum p_{\alpha} l_{\alpha} = 0. 
\]

Next we show that each $\mu$ for which $\vert \mu \vert > 1$ is real and positive. Let $l_{\alpha}$ satisfy (\ref{E:letterE}) with  $\vert \mu \vert > 1$; then the eigenvalue for $\overline{l}_{\alpha}$ is $\overline{\mu}$ and so 
\[
	\sum_\beta{} \left( \mathcal{Q}^{r} \right)_{\alpha \beta}  \left( 1 +\varepsilon \left( l_{\beta} + \overline{l}_{\beta} \right) \right) = 1 + 2 \varepsilon \; \Re \; \mu^{r}l_{\alpha}.
\]
If $\varepsilon > 0$ has been chosen so small that $\Re \, \varepsilon l_{\beta} > -1/2$ for all $\beta$ then the L.H.S. is positive for the coefficients in the matrix are positive, whereas the R.H.S. is negative for suitably chosen $\mu$, unless $l_{\alpha} = 0$. If now $\mu > 1$ we may take it that $l_{\alpha}$ is real for each $\alpha$. As it must satisfy $\sum p_{\alpha} l_{\alpha} = 0$ it is negative for some $\alpha$, but then
\[
	\sum_{\beta} \left(\mathcal{Q}^{r} \right)_{\alpha \beta} \left( 1 + \varepsilon \; l_{\beta} \right) = 1 + \varepsilon \mu^{r} l_{\alpha},
\]
%
%
and if $\varepsilon$ is chosen so that $1 + \varepsilon \; l_{\beta} > 0$ for all $\beta$ the L.H.S is positive whereas the R.H.S is negative for sufficiently large $r$. 

All the eigenvalues therefore satisfy $\vert \mu \vert \leq 1$ as the eigenvalues are all different in modulus this means that $\vert \mu \vert < 1$ except for one value of $\mu$. Then as $r \to \infty, \mathcal{M}^{r}$ tends to a matrix which has only one element different from 0, and that a 1 on the diagonal, say in position $\sigma\sigma$. 

Calling this matrix $\mathcal{X}_{\sigma}$ the series of matrices $\mathcal{Q}^{r}$ tends to the matrix $\mathcal{U}^{-1}\mathcal{X_{\sigma}U}$. This matrix is the one and only one $\mathcal{Y}$ which satisfies $\mathcal{YQ = Y, Y^{\text{2}} = Y, Y \ne \text{0}}$ and is therefore the one whose $\alpha \beta$ coefficient is $\mathcal{P}_{\beta}$.

%
%
\subsection{Another probability problem}
There is another probability problem that arises in connection with simple transpositions. With a message of length $L$, and a key length of $K$ what is the probability that the $m$th letter will be at the bottom of a column? Let $D$ be the length of the short columns \textit{i.e.} $D = \left[ L/K \right]$, and let $E = L - DK$. Then if the $m$th letter is at the bottom of the $w$th column we must have
\[
	\frac{m}{D + 1} \leq w  \leq \frac{m}{D},
\]
and there will be $(D + 1)w- m$ short and $m - Dw$ long columns among these first $w$ columns. There are\footnote{ Turing is using Binomial Coefficient notation in this section;
\[
	\binom{n}{k} =  C(n, k) = \frac{P(n,k)}{P(k, k)} = \frac{n!}{(n-k)! k!}
\]
}

\[
	\binom{w}{m - Dw}  \binom{K -w}{E - m + Dw}
\]
ways in which the short and long columns can be arranged consistently with this, and altogether $\binom{K}{E}$ ways in which the columns can be arranged, so that the probability of the $m$ the letter being at the bottom of a column is
\large{
\[
	\sum \limits_{\left( m/D+1 \right) \le w \le \left( m/D \right)} \binom{w}{m - Dw} \binom{K - w}{E - m + Dw} \Bigg/ \binom{K}{E}.
\]}
\normalsize
There will normally be very few terms in the sum. 

Let us take the case of the message of length 133 and consider the 45th letter, assuming the key length is between 10 and 20 (inclusive). $L_{O}{}^{B}. \quad L = 133, m =~45.$

\begin{align*}
&K = 10, \quad D = 13, \quad E = 3, \quad \frac{m}{D + 1} = 3+, \quad \frac{m}{D}  = 3+ &&\textit{no terms}\\
&K = 11, \quad D = 12, \quad E = 1, \quad \frac{m}{D + 1} = 3+,  \quad \frac{m}{D}  = 3+ &&\textit{no terms}\\
&K = 12, \quad D = 11, \quad E = 1, \quad \frac{m}{D + 1} = 3+, \quad \frac{m}{D}  = 4+ && \\
& \quad  \textit{only terms $w = 4$, $m-Dw = 1$  \qquad \qquad  probability is:} && \binom{4}{1}  \binom{8}{0} \Bigg/ \binom{12}{1} = \frac{4}{12} \\
&K = 13, \quad D = 10, \quad E = 3, \quad \frac{m}{D + 1} = 4+, \quad \frac{m}{D}  = 4+ &&\textit{no terms}\\\\
&K = 14, \quad D = 9, \quad E = 7, \quad \frac{m}{D + 1} = 4+, \quad \frac{m}{D}  = 5 \quad && \\
& \quad \textit{only terms $w = 5$,  $m-Dw = 0$ \qquad \qquad probability is:}  &&\binom{5}{0}  \binom{9}{7} \Bigg/  \binom{14}{7}  = \frac{3}{286} = 0.0105,  \\
&K = 15, \quad D = 8, \quad E = 13, \quad \frac{m}{D + 1} = 5, \quad \frac{m}{D}  = 5+ \\
&\quad \textit{only terms $w = 5$, $m-Dw = 5 $  \qquad \qquad probability is:} &&\binom{5}{5}  \binom{10}{8} \Bigg/  \binom{15}{13}  = \frac{3}{7} = 0.428,  \\
&K = 16, \quad D = 8, \quad E = 5, \quad \frac{m}{D + 1} = 5+, \quad \frac{m}{D}  = 5+ \\
&\quad \textit{only terms $w = 5$, $m-Dw = 5$ \qquad \qquad probability is:} &&\binom{5}{5}   \binom{11}{0}  \Bigg/  \binom{16}{5} = \frac{1}{4368} = 0.000229,  \\
&K = 17, \quad D = 7, \quad E = 14, \quad \frac{m}{D + 1} = 5+, \quad \frac{m}{D}  = 6+ \\
&\quad \textit{only terms $w = 6$, $m-Dw = 3$ \qquad \qquad probability is:}  &&\binom{6}{3}  \binom{11}{11} \Bigg/  \binom{17}{14}  = \frac{1}{34} = 0.0307 \\
& && \qquad \qquad \qquad {(Editor -   1/34 \textnormal{ is } 0.0294.)} \\
&K = 18, \quad D = 7, \quad E = 7, \quad \frac{m}{D + 1} = 5+, \quad \frac{m}{D}  = 6+ \\
&\quad \textit{only terms $w = 6$, $m-Dw = 3$ \qquad \qquad probability is:}  &&\binom{6}{3}  \binom{12}{4}  \Bigg/  \binom{18}{7}  = \frac{4950}{15912} = 0.311,  \\
&K = 19, \quad D = 7, \quad E = 0, \quad  \qquad \qquad \qquad \textit{ probability is:} && = 0\\
&K = 20, \quad D = 6, \quad E = 13, \quad \frac{m}{D + 1} = 6+, \quad \frac{m}{D}  = 7+ \\
&\quad \textit{only terms $w = 7$, $m-Dw = 3$ \qquad \qquad probability is:} &&\binom{7}{3}  \binom{13}{4}  \Bigg/  \binom{20}{7}  = \frac{35 \times143}{15504} = 0.323.
\end{align*}

\vspace{10pt}
\hrule
\vspace{4pt}
\begin{center}
	$\infty$
\end{center}
\vspace{4pt}
\hrule

\end{document}